\documentclass[12pt]{article}

\usepackage[russian,ukrainian,english]{babel}

\usepackage{amssymb}
\usepackage{graphics}

\newcommand{\qed}{$\;\;\;\Box$}
\newenvironment{proof}{\par\smallbreak{\sl\bf Proof.~}}
{\unskip\nobreak\hfill \qed \par\medbreak}

\newcounter{claim}
\renewcommand{\theclaim}{\arabic{claim}}
{\par\medskip\par}

{\qed\par\smallbreak}
\newcommand{\N}{{\mathbb N}}
\newcommand{\R}{{\mathbb R}}

\newcommand{\LL}{{\cal L}}

\newcommand{\beq}{\begin{equation}}
\newcommand{\ee}{\end{equation}}

\renewcommand{\d}{\partial}

\newtheorem{thm}{Theorem}[section]
\newtheorem{lem}[thm]{Lemma}
\newtheorem{defn}[thm]{Definition}
\newtheorem{cor}[thm]{Corollary}

\newtheorem{ex}[thm]{Example}

\newcommand{\al}{\alpha}
\newcommand{\be}{\beta}
\newcommand{\ga}{\gamma}

\newcommand{\eps}{\varepsilon}
\newcommand{\vphi}{\varphi}
\newcommand{\la}{\lambda}
\newcommand{\om}{\omega}
\newcommand{\io}{\iota}

\newcommand{\reff}[1]{(\ref{#1})}      

\newcommand{\diag}{\mathop{\rm diag}\nolimits}

\title{
Exponential  stability of solutions to perturbed superstable wave equations}

\newcounter{thesame}
\author{
I.~Kmit
\thanks{Institute of Mathematics, Humboldt University of Berlin. On leave from the
Institute for Applied Problems of Mechanics and Mathematics,
Ukrainian National Academy of Sciences. {\small   E-mail:
{\tt kmit@mathematik.hu-berlin.de}}}
\ \ \ N.~Lyul'ko \thanks{Sobolev Institute of Mathematics, Russian Academy of Sciences and
Novosibirsk State University, Russia.
{\small   E-mail:
{\tt natlyl@mail.ru}}
}}
\date{}

\begin{document}

\maketitle

\begin{abstract}
\noindent 
The paper deals with  initial-boundary value problems
 for the linear wave equation whose solutions  stabilize to zero in a finite time.
We prove that problems in this class remain 
 exponentially stable in $L^2$ as well as in $C^2$ under  small bounded perturbations 
of the wave operator. 
To show this for $C^2$, we prove a smoothing result implying that
the solutions to the perturbed problems
 become eventually $C^2$-smooth for any $H^1\times L^2$-initial data.
\end{abstract}



\section{Motivation and our results}\label{sec_1}

 A linear system
\begin{equation}\label{ss0}
\frac{d}{dt}x(t)=A(t)x(t), \quad x(t)\in X \quad (0\le t \le \infty),
\end{equation}
on a Banach space $X$
is called  \textit{exponentially stable} if there exist positive reals
$\gamma$ and $M=M(\gamma)$ such that every solution $x(t)$ 
satisfies the estimate
\begin{equation}
\|x(t)\|\le M e^{-\gamma t} \|x(0)\|, \quad t\ge 0,\label{cor1}
\end{equation}
where  $\|\cdot\|$ denotes the norm in $X$. 

The papers  \cite{bal105,cr13} address a stronger property of exponentially stable systems, known as superstability. They consider the Cauchy problem for the  {\it autonomous} version of (\ref{ss0}). Moreover,
$A: X \rightarrow X$ is supposed to be  the infinitesimal generator of a strongly continuous semigroup  $T(t)$; see \cite{hp57,Paz}. 
A semigroup $T(t)$ is called {\it superstable}  \cite{bal105} 
if its stability index 
is $-\infty$, that is
$$
\lim_{t\to \infty}\frac{\log\|T(t)\|}{t}=-\infty.
$$
The superstability property 
implies that the system is exponentially stable and, moreover, the estimate (\ref{cor1})
holds for every $\gamma> 0$. 

The notion of superstability can be appropriately extended to second-order equations, studied in this paper.

An important subclass of  superstable systems  consists of the systems
whose solutions  stabilize to zero 
after some  time. The time of the stabilization is called  a {\it finite time extinction}. 
The following simple example shows which type of boundary conditions cause the superstability property for the wave 
equation in the one-dimensional case, see  \cite{bal105,perroll}:  
\begin{equation}\label{rr1}
w_{tt}  -a^2w_{xx}  = 0, \quad (x,t)\in \Pi, 
\end{equation}
\begin{equation}\label{rr2}
\begin{array}{ll}
w(0,t)=0, \quad t\in(0,\infty),\\
 (w_t  +a w_x)(1,t)  = 0, \quad t\in(0,\infty),
\end{array}
\end{equation}
\begin{equation}\label{rr3}
\begin{array}{ll}
w(x,0)=w_0(x),\quad x \in[0,1],\\ w_t(x,0)=w_1(x),\quad x \in[0,1],
\end{array}
\end{equation}
where  $a$ is arbitrary positive constant and  $\Pi=(0,1)\times(0,\infty)$. Constructing  solutions  by the method of characteristics shows that all solutions stabilize to zero for
 $t>2/a$. 
A class of  boundary conditions ensuring the superstability property for the wave operator in  the multidimensional case is described in
 \cite{Majda}.
Such boundary conditions naturally appear,
for example, in the scattering theory \cite{lax1,lax2} and in the  control theory \cite{moulay}.

A comprehensive review of the available results on asymptotic
 behavior of solutions to linear and quasi-linear first order hyperbolic problems can be 
found in  \cite{bastin,KL}. 
The present paper concerns superstable  initial-boundary value problems for the one-dimensional wave equation.
Specifically, 
we consider the problem 
(\ref{rr1}), (\ref{rr3}) with the boundary conditions either 
\begin{equation}\label{rr4}
\begin{array}{ll}
\displaystyle
w(0,t)=p(w_t +a w_x)(0,t),\quad  t\in(0,\infty),
\nonumber\\
\displaystyle
(w_t +a w_x)(1,t)=0,\quad  t\in(0,\infty),
\nonumber
\end{array}
\end{equation}
or
\begin{equation}\label{rr41}
\begin{array}{ll}
\displaystyle
w(1,t)=q(w_t -a w_x)(1,t), 
\quad  t\in(0,\infty),
\nonumber\\
\displaystyle
(w_t -a w_x)(0,t)=0,\quad  t\in(0,\infty),
\nonumber
\end{array}
\end{equation}
where $p$ and $q$ are arbitrary constants.
Note that the finite time extinction for the problems under consideration equals $2/a$ (see Section 3.1). In this paper we will investigate the problem (\ref{rr1}), (\ref{rr4}), (\ref{rr3})
(the same result is true for the problem (\ref{rr1}), (\ref{rr41}), (\ref{rr3})
as well).

Along with the equation  (\ref{rr1}) we will consider its perturbed version, namely
\begin{equation}\label{rr6}
w_{tt}-a^2w_{xx} +c(x,t)w = 0, \quad (x,t)\in \Pi,  
\end{equation}
where $c$ is a two times continuously differentiable function such that $c$ itself and its first order and second order derivatives 
are bounded on $\overline\Pi$.  

Given  $T>0$, write $\Pi_T=(0,1)\times(0,T)$.  By $H^1(0,1)$ we  denote the space of  functions 
$u: (0,1)\to R$ such that $u\in L^2(0,1)$ and its weak derivative  $u^\prime\in L^2(0,1)$. 
From \cite{KL} it follows  that for any initial functions $w_0\in H^1(0,1)$ and $ w_1\in L^2(0,1)$
the problem  (\ref{rr6}), (\ref{rr4}), (\ref{rr3}) has a unique $L^2$-generalized solution $w\in C([0,\infty),L^2(0,1))$
(see Section 2 for the definition). For this solution we prove the following perturbation theorem.

\begin{thm}\label{th1} 
Let $w_0\in H^1(0,1)$ and $ w_1\in L^2(0,1)$.
Then  for  any $\gamma >0$ there exist $\varepsilon >0$ and $M=M(\gamma)$ such that,
whenever  $\sup_{(x,t)\in \Pi}|c|<\varepsilon$, the $L^2$-generalized solution $w(x,t)$ to the problem
(\ref{rr6}),  (\ref{rr4}), (\ref{rr3})   fulfills the bound
\begin{equation}\label{bound1}
\|w(\cdot,t)\|_{L^2(0,1)}\le 
M e^{-\gamma t}\max(\|w_0\|_{H^1(0,1)},\|w_1\|_{L^2(0,1)}), \quad t>0.
\end{equation}
\end{thm}

This means that sufficiently small bounded perturbations of the zero order part of the wave equation (\ref{rr1})
 presurve the exponential stability in $L^2(0,1)$ of the unperturbed problem (\ref{rr1}),  (\ref{rr4}), (\ref{rr3}).

\begin{defn}\label{def}  
The problem (\ref{rr6}), (\ref{rr4}), (\ref{rr3}) is called {\rm  smoothing} from $L^2(0,1)$ to $C^2([0,1])$ if there is
 $T>0$ such that, given 
   $w_0\in H^1(0,1)$ and $w_1\in L^2(0,1)$, the $L^2$-generalized solution $w(x,t)$ to this problem belongs to  $C^2([0,1])$
 whenever $t\ge T$.  
\end{defn}

We are prepared to formulate the following smoothing result.

\begin{thm}\label{th2} 
The  problem (\ref{rr6}), (\ref{rr4}), (\ref{rr3}) is smoothing from $L^2(0,1)$ to $C^2([0,1])$ with $T=6/a$. 
Moreover, given $w_0\in H^1(0,1)$ and $w_1\in L^2(0,1)$, the $L^2$-generalized solution $w(x,t)$  fulfills the bound
\begin{equation}\label{bound2}
\|\partial^{\alpha,\beta}_{x,t} w(\cdot,t)\|_{C([0,1])}\le
M_1 e^{\omega t}\max(\|w_0\|_{H^1(0,1)},\|w_1\|_{L^2(0,1)}),\qquad t\ge T,
\end{equation}
where $\alpha+\beta\le 2$ and $M_1, \omega$ are constants not depending on $t$, $w_0$, and $w_1$.
\end{thm}

\begin{cor}\label{cor} 
Let $w_0\in H^1(0,1)$ and $ w_1\in L^2(0,1)$.
Then
 for  any $\gamma >0$ there exist $\varepsilon >0$ and $M_2=M_2(\gamma)\ge M$ such that,
whenever $\sup_{T>0}(\|c\|_{C^2(\overline{\Pi_T})})<\varepsilon$,
the $L^2$-generalized solution $w(x,t)$ to the problem (\ref{rr6}),  (\ref{rr4}), (\ref{rr3})  fulfills the bound
\begin{equation}\label{bound3}
\|\partial^{\alpha,\beta}_{x,t} w(\cdot,t)\|_{C([0,1])}\le 
M_2 e^{-\gamma t}\max(\|w_0\|_{H^1(0,1)},\|w_1\|_{L^2(0,1)}), \qquad t>T,
\end{equation}
where $T=6/a$,  $\alpha+\beta\le 2$, and $M_2, \gamma$ are constants not depending on $t$, $w_0$, and $w_1$.
\end{cor}

\section{Reduction to first order systems and $L^2$-generalized solutions}\label{sec_2} 

For proving Theorems \ref{th1} and \ref{th2}, we will apply our results from \cite{KL}, where we study the
exponential stability of the solutions to the initial-boundary value problems for the general first order nonautonomous 
hyperbolic systems, that are small bounded perturbations of the corresponding superstable systems. To this end, let us
introduce new variable $u=\partial_t w+a\partial_x w$ and
rewrite the perturbed problem as the following $2\times 2$-first order system:
\begin{equation}\label{ss1}
\partial_t w+a\partial_x w-u=0, \qquad \partial_t u-a\partial_x u+c(x,t)w=0,
\end{equation}
\begin{equation}\label{ss2}
 w(0,t)=pu(0,t), \qquad  u(1,t)=0,
\end{equation}
\begin{equation}\label{ss21}
 w(x,0)=\varphi_1(x),\quad u(x,0)=\varphi_2(x),
\end{equation}
where
$$
\varphi_1(x)=w_0(x), \quad \varphi_2(x)=w_1(x)+a\partial_x w_0(x).
$$
It is evident that the $C^2$-function $w$ is a classical solution to the problem  (\ref{rr6}), (\ref{rr4}), (\ref{rr3})
if and only if the $C^2$-vector function  $(w, u)=(w,\partial_t w+a\partial_x w)$ is 
a classical solution to the problem  (\ref{ss1})--(\ref{ss21}).

As it follows from \cite{ijdsde}, for $w_0\in C^3[0,1]$ and $w_1\in C^2[0,1]$ satisfying the zero order, the first order 
and the second order compatibility conditions between (\ref{ss2}) and  (\ref{ss21}) the problem 
(\ref{ss1})--(\ref{ss21}) has a unique classical solution $(w,u)\in C^2([0,1])\times C^2([0,1])$, where the first component $w$ is a classical solution to the problem  (\ref{rr6}), (\ref{rr4}), (\ref{rr3}). Note that then the 
 compatibility conditions up to the second order  between (\ref{rr4}) and  (\ref{rr3}), see e.g.  \cite{Lyl}, follow immediately from the
equality $u=\partial_t w+a\partial_x w$.

Let us introduce the notion of an $L^2$-generalized solution. Fix  arbitrary  $w_0\in H^1(0,1)$ and $ w_1\in L^2(0,1)$ and 
sequences  $\varphi_1^l\in C_0^\infty([0,1])^n$, $\varphi_2^l\in C_0^\infty([0,1])^n$
such that   $\varphi_i^l\to\varphi_i$ in $L^2(0,1)$, $i=1, 2$. Note that, due to the fact that  $\varphi_i^l$ are compactly supported for all 
$l\in N$, they satisfy 
the  compatibility conditions up to the second order between (\ref{ss2}) and  (\ref{ss21}). By \cite{ijdsde},
given  $l\in N$, the problem (\ref{ss1})--(\ref{ss21}) has a unique classical solution, say $(w^l,u^l)$, belonging to
 $C^2([0,1])\times C^2([0,1])$.
In \cite{KL} it is proved that these solutions satisfy the following estimate:
\begin{equation}\label{bound40}
\max(\|w^l(\cdot,t)\|_{L^2(0,1)}, \|u^l(\cdot,t)\|_{L^2(0,1)})\le 
M_3 e^{A t}\max_{i=1,2}(\|\varphi_i\|_{L^2(0,1)}),\quad \mbox{ for all } t>0,
\end{equation}
for some constants $M_3, A$, not depending on $t$, $\varphi_1$,  $\varphi_2$, and $l\in N$. Hence, there exist
unique functions  $w, u \in C([0,\infty), L^2(0,1))$ such that 
$$
\|w(\cdot,\theta)-w^{l}(\cdot,\theta)\|_{L^2(0,1)}\to 0, \,\, \|u(\cdot,\theta)-u^{l}(\cdot,\theta)\|_{L^2(0,1)}\to 0\quad 
\mbox{ as } l\to\infty,
$$
uniformly in $\theta$ varying in the range $0\le\theta\le t$, for every $t>0$. 
The vector-function ($w, u$) is called an  $L^2$-{\it generalized  solution} 
to the problem (\ref{ss1})--(\ref{ss21}),  while the function $w$
is called an   $L^2$-{\it generalized  solution} to the problem (\ref{rr6}), (\ref{rr4}), (\ref{rr3}).
Furthermore, the following estimate is true:
\begin{equation}\label{bound4}
\|w(\cdot,t)\|_{{L^2(0,1)}}\le 
M_3 e^{A t}\max(\| w_0\|_{H^1(0,1)},\|w_1\|_{L^2(0,1)})\quad  \mbox{ for all } t>0,
\end{equation}
what follows from (\ref{bound40}). 

\section{Proofs of the main results}

\subsection{Proof of Theorem \ref{th1} }

Consider the problem (\ref{ss1})--(\ref{ss21}). Due to the method of characteristics, the classical solution to this problem 
fulfills the following system of integral equations:
\begin{equation}\label{ss14} 
w(x,t)= \int _{t-x/a}^t u(a(\tau-t)+x,\tau)\,d\tau+p u\left(0,t-x/a\right), \quad t>\frac{x}{a},
\end{equation}
\begin{equation}\label{ss15} 
u(x,t)= 
-\int _{t+(x-1)/a}^t [cw](x+a(t-\tau),\tau)\,d\tau,\,\quad t+\frac{x}{a}>\frac{1}{a}.
\end{equation}
Note that if $c=0$ then  the finite time extinction for this problem equals $2/a$.

Letting  $t>4/a$ and substituting (\ref{ss15}) into (\ref{ss14}), we get an integral equation for $w$, namely
$$
w(x,t)=- \int _{t-x/a}^t \int _{2\tau-t+(x-1)/a}^\tau\,[cw](a(\tau-\xi)+x,\xi)\,d\xi d\tau-p \int _{t-(x+1)/a}^{t-x/a} [cw](a(t-\tau)-x,\tau)\,d\tau.
$$
Set $W(t)=\|w(\cdot,t)\|_{L^2(0,1)}$. Then the last equality entails
\begin{equation}\label{ss16} 
W(t)\le K\sup_{(x,t)\in \Pi}|c(x,t)|  \int _{t-\frac{4}{a}}^t W(\tau)\,d\tau, \quad t>\frac{4}{a},
\end{equation}
the constant $K$ being independent on $t, w$ but only on   $p$. On the other side, the estimate (\ref{bound4}) causes the bound 
\begin{equation}\label{bound42}
W(t)\le 
K_1\max(\| w_0\|_{H^1(0,1)},\|w_1\|_{L^2(0,1)}),\quad    0\le t \le \frac{4}{a},
\end{equation}
with the constant $K_1=M_3 e^{\frac{4A }{a}}$. On the account of
(\ref{ss16}), (\ref{bound42}) and   \cite[Lemma 5.1]{KL},   for  any $\gamma>0$ there exist $\varepsilon>0$ and $M=M(\gamma)$ such that,
whenever  $\sup_{(x,t)\in \Pi}|c(x,t)|<\varepsilon$, for the classical solution $w$ to the perturbed problem  (\ref{rr6}), (\ref{rr4}), (\ref{rr3})  the following bound is true:
$$
W(t)\le 
M e^{-\gamma t}\max(\| w_0\|_{H^1(0,1)},\|w_1\|_{L^2(0,1)}), \quad t>0. 
$$
Finally, the statement of the theorem follows from the definition of an  $L^2$-{\it generalized  solution} 
to the problem (\ref{rr6}), (\ref{rr4}), (\ref{rr3}).

\subsection{Proof of Theorem \ref{th2}}

First consider the problem (\ref{ss1})--(\ref{ss21}) where the system (\ref{ss1}) is replaced by the corresponding
decoupled system, namely
\begin{equation}\label{ss11}
\begin{array}{ll}
\partial_t w+a\partial_x w=0, \quad \partial_t u-a\partial_x u=0,\\
 w(0,t)=pu(0,t), \quad  u(1,t)=0,\\
 w(x,0)=\varphi_1(x), \quad u(x,0)=\varphi_2(x).
\end{array}
\end{equation}
This system is easily checked to be superstable. Indeed, for the classical solution  $(w,u)$ to (\ref{ss11}) the method of characteristics
gives the formulas
$$
w(x,t) =w(0, t-\frac{x}{a})=pu\left(0,t-\frac{x}{a}\right), \quad t>\frac{x}{a},
$$
$$
u(x,t)= u\left(1,t+\frac{x-1}{a}\right), \quad  t+\frac{x}{a}>\frac{1}{a}.
$$
Since $u(1,t)=0$, we conclude that the  finite time extinction for this problem is $2/a$. 

Then the smoothing property of the problem (\ref{ss1})--(\ref{ss21}) 
immeadiately follows from \cite[Theorem 12]{Km} and \cite[Theorem 2.5]{KL}. More specifically, 
there exists
a positive real $T$ such that for any $\varphi_1, \varphi_2 \in L^2(0,1)$ and any $C^2$-function  $c$  the $L^2$-generalized solution $(w,u)$   to the problem (\ref{ss1})--(\ref{ss21}) is two times continuously differentiable for $t>T$. Moreover, it fulfills the following estimate:
$$
\max(\|\partial^{\alpha,\beta}_{x,t} w(\cdot,t)\|_{C[0,1]}, \|\partial^{\alpha,\beta}_{x,t}u(\cdot,t)\|_{C([0,1])})\le
M_1 e^{A t}\max_{i=1,2}(\|\varphi_i\|_{L^2(0,1)}), \quad t>T,
$$
where  $\alpha+\beta\le 2$  and the constants $M_1, A$ are independent of $t$ and $\varphi_1$, $\varphi_2$. This entails that
the first component $w$ of this solution, which is the $L^2$-generalized solution to the original problem 
 (\ref{rr6}), (\ref{rr4}), (\ref{rr3}), fulfills Theorem 3, as desired.

Finally,  to prove Corollary \ref{cor}, we follow the argument used to prove \cite[Theorem 2.7($ii$)]{KL}. The estimate
(\ref{bound3}) will then follow from the estimates (\ref{bound1}) and (\ref{bound2}).

\subsection{Generalization}

Our results given by Theorems\ref{th1} and \ref{th2} can be easily extended to the second order equations 
involving  first order terms, of the following type:
\begin{equation}\label{nl1}
(\partial_t -a\partial_x+a_1(x,t))(\partial_t +a\partial_x)w +c(x,t)w=0,
\end{equation}
where  the coefficient $a_1$ is a two times continuously differentiable function such that $a_1$ itself and its first order and second order derivatives are bounded on $\overline\Pi$. If the coefficient $a$ is  positive, then
 in the domain $\Pi$ we consider  the initial-boundary value problem (\ref{nl1}), (\ref{rr4}), (\ref{rr3});
if $a$ is negative, then we consider the problem (\ref{nl1}), (\ref{rr41}), (\ref{rr3}).

For the problem (\ref{nl1}), (\ref{rr4}), (\ref{rr3}), in the proofs we encounter the following minor changes.
Similarly to the above, the non-perturbed problem (\ref{nl1}), (\ref{rr4}), (\ref{rr3}) (with $c\equiv 0$) has the finite time 
 extinction $2/a$, while the first order system  corresponding to (\ref{nl1}) now reads  as
\begin{equation}\label{ss111}
\partial_t w+a\partial_x w=u, \qquad \partial_t u-a\partial_x u+a_1(x,t) u+c(x,t)w=0
\end{equation}
and the boundary and the initial conditions  (\ref{rr4}) and (\ref{rr3}) read as (\ref{ss2}), (\ref{ss21}).
The solution to the first order problem (\ref {ss111}), (\ref{ss2}), (\ref{ss21}) satisfies the following integral system:

$$
w(x,t)= \int _{t-\frac{x}{a}}^t u(a(\tau-t)+x,\tau)\,d\tau+p u(0,t-\frac{x}{a}), \quad t>\frac{x}{a},
$$
$$
u(x,t)= e^{\int _1^x   \frac{a_1(\xi,t+\frac{x-\xi}{a})}{a}\,d\xi}
\int _{t+\frac{x-1}{a}}^t [-cw e^{\int_\eta^1 \frac{a_1(\xi,\tau+\frac{\eta-\xi}{a})\,d\xi}{a}}](\eta,\tau)|_{\eta=x+a(t-\tau)}d\tau,\,\quad t+\frac{x}{a}>\frac{1}{a}.
$$

\section*{Acknowledgments}
 Irina Kmit was supported by the 
VolkswagenStiftung Project ``Modeling, Analysis, and Approximation Theory toward Applications in 
Tomography and Inverse Problems''. Natalya Lyulko was supported by the Presidium of the Russian Academy of Sciences 
(basic research program, number I.5P).


\begin{thebibliography}{10}

\bibitem{bastin}G. Bastin,  J.-M. Coron, {\it Stability and Boundary Stabilization of 1-D Hyperbolic Systems},
Progress in Nonlinear Differential Equations and Their Applications {\bf 88}, 
Birkh\"auser, 2016.


\bibitem{bal105} A.V. Balakrishnan, Superstability of systems,
{\it  Applied Mathematics and Computation} {\bf 164(2)} (2005), 321--326.

\bibitem{cr13}  D. Creutz, M. Mazo Jr., C. Preda, Superstability and finite time  extinction  
for $C_0$-semigroups, (2013). {\it  E-print}: https://arxiv.org/abs/0907.4812.

\bibitem{hp57} E. Hille, R. Phillips,  {\it Functional  analysis and semi-groups},
Providence, 1957.

\bibitem{ijdsde}
I. Kmit,  Classical solvability of nonlinear initial-boundary problems for
first-order hyperbolic systems, {\it Intern. J. Dynamic Syst. Differ. Equat.} {\bf 1(3)} (2008), 191--195.

\bibitem{Km} I. Kmit, Smoothing solutions to initial-boundary problems  for first-order
hyperbolic systems, {\it Applicable Analysis} {\bf 90(11)} (2011), 1609--1634.


\bibitem{KL} I. Kmit, N. Lyulko, Perturbations of superstable linear hyperbolic systems, (2017).
{\it E-print:} https://arxiv.org/abs/1605.04703

\bibitem{lax1} P. Lax, R. Phillips,
{\it Scattering theory}, Academic Press, New York, 1967.

\bibitem{lax2} P. Lax, R. Phillips, Scattering theory for
the dissipative wave equation, {\it Indiana Univ. Math. J.}, {\bf 24} 1975,
1119--1138.

\bibitem{Lyl} N. Lyulko, Increasing
smoothness of solutions to mixed problem for the wave equation on the plane, {\it J. Math. Phys., Anal., Geom.}
{\bf 11(2)}  (2004), 169--176.

\bibitem{Majda}  A. Majda,
Disappearing solutions for the dissipative wave equation,
{\it Indiana Univ. Math. J.} {\bf 24} (1975), 1119--1133.

\bibitem{moulay}E. Moulay, W. Perruquetti,
Finite time stability and stabilization:  state of art, {\it Advances in variable structure
and sliding mode control}, 23–41, Lecture Notes in Control and Inform. Sci., {\bf 334}  2006, Springer, Berlin,  23-–41.

\bibitem{Paz} 
A. Pazy,
{\it Semigroups of operators and applications to partial differential equations},
Springer-Verlag, Berlin,  1983. 

\bibitem{perroll} V. Perrollaz, L. Rosier, Finite-Time Stabilization of $2\times2$ Hyperbolic Systems on Tree-Shaped Networks, {\it  SIAM Journal on Control and Optimization} {\bf 52(1)}
(2014), 143--163.


\end{thebibliography}
\end{document}